\documentclass[11pt,a4paper]{article}% For Computer Modern Fonts, or,
\usepackage{latexsym,amssymb, amsthm}
\usepackage[centertags]{amsmath}
\usepackage{epsfig,pifont}
\usepackage[inactive]{srcltx}  %SRC Specials for DVI search

\usepackage{graphics,helvet,times,mathptm,epic,eepic}

\usepackage{color}

\usepackage{newlfont}

\usepackage{amsfonts,bbm}

\newenvironment{keywords}{ \noindent {\small\bf Key Words}:}{ }

\def\bd{\begin{description}}
\def\ed{\end{description}}

\def\beq{\begin{equation}}
\def\eeq{\end{equation}}
\def\bea{\begin{eqnarray}}
\def\eea{\end{eqnarray}}
\def\beas{\begin{eqnarray*}}
\def\eeas{\end{eqnarray*}}

\newtheorem{theorem}{Theorem}[section]

\newtheorem{corollary}{Corollary}

\theoremstyle{remark}
\newtheorem{example}{Example}[section]

\begin{document}
%\begin{article}

\title{\textbf{\textsc{Numerical computations and mathematical modelling with  infinite and infinitesimal numbers}}}

%% Author name
\newcommand{\nms}{\normalsize}
\author{  {   \bf Yaroslav D. Sergeyev\footnote{Yaroslav D. Sergeyev, Ph.D., is Distinguished Professor  at the University of Calabria, Rende, Italy.
 He is also Full Professor (a part-time contract) at the N.I.~Lobatchevsky State University,
  Nizhni Novgorod, Russia and Affiliated Researcher at the Institute of High Performance
  Computing and Networking of the National Research Council of Italy.}
 % $\hspace{2mm}^,$
 % \footnote{The author thanks the anonymous referees for their useful suggestions.}
 }\\ \\ [-2pt]
      \nms Dipartimento di Elettronica, Informatica e Sistemistica,\\[-4pt]
       \nms   Universit\`a della Calabria,\\[-4pt]
       \nms 87030 Rende (CS)  -- Italy\\ \\[-4pt]
       \nms http://wwwinfo.deis.unical.it/$\sim$yaro\\[-4pt]
         \nms {\tt  yaro@si.deis.unical.it }
}

\date{}

\maketitle

\begin{abstract}
Traditional computers work with finite numbers. Situations where
the usage of infinite or infinitesimal quantities is required are
studied mainly theoretically. In this paper, a recently introduced
computational methodology (that is not related to the non-standard
analysis) is used to work with finite, infinite, and infinitesimal
numbers \textit{numerically}. This can be done on  a new kind of a
computer -- the Infinity Computer -- able to work  with all these
types of numbers.
  The new computational tools both give
possibilities to execute computations of a new type and open new
horizons for creating new mathematical models where a
computational usage of infinite and/or infinitesimal numbers can
be useful. A number of numerical examples showing the potential of
the new approach and  dealing with divergent series, limits,
probability theory, linear algebra, and calculation of volumes of
objects consisting of parts of different dimensions are given.
 \end{abstract}


%\newpage

\section{Introduction}
\label{s0}

The point of view on infinity accepted nowadays takes its origins
from the famous ideas of Georg Cantor (see \cite{Cantor}).
Different   generalizations of the traditional arithmetic for
finite numbers to the case of infinite and infinitesimal numbers
have been proposed in literature (see, e.g.,
\cite{Benci,Cantor,Conway,Loeb,Robinson} and references given
therein). However, these generalizations are quite different with
respect to the finite arithmetic we are used to deal with.
Moreover, very often they leave undetermined many operations where
infinite numbers take part (e.g., $\infty-\infty$,
$\frac{\infty}{\infty}$, etc.) or use representation of infinite
numbers based on infinite sequences of finite numbers.

In spite of these crucial difficulties and due to the enormous
importance of the concept of infinite and infinitesimal in
Science, people try to introduce these notions in their work with
computers. Thus, the IEEE Standard for Binary Floating-Point
Arithmetic (IEEE 754) being the most widely-used standard for
floating-point computation  defines formats for representing
  special values for positive and
negative infinities and NaN (Not a Number)  (see also
incorporation of these notions   in the interval analysis
implementations   e.g., in \cite{Walster_1}). The IEEE infinity
values   can be the result of arithmetic overflow, division by
zero, or other exceptional operations. In   turn, NaN  is a value
or symbol that can be produced as the result of a number of
operations including that involving zero, NaN itself, and
infinities.

Recently, a new applied point of view on infinite and
infinitesimal numbers has been introduced in
\cite{Sergeyev,Poland,chaos,spirals,informatica}. With respect to
the IEEE 754 standard, the new approach significantly extends the
variety of operations that can be done with infinity. It gives a
possibility to work with various infinite  and infinitesimal
quantities \textit{numerically} by using a new kind of a computer
-- the Infinity Computer -- introduced in
\cite{Sergeyev_patent,www,Poland}. A number of applications
related to the usage of the new approach for studying fractals
(being one of the main scientific interests of the author (see,
e.g., \cite{Strongin_Sergeyev,Strongin&Sergeyev(2003)}) has been
discovered (see \cite{chaos,spirals}).

The new approach is not related to non-standard analysis ideas
from \cite{Robinson} and does not use Cantor's ideas (see
\cite{Cantor}) either. The Infinity Computer works with infinite
and infinitesimal numbers numerically using the following
methodological principles having a strong applied character (see
survey \cite{informatica} for a detailed discussion on the new
approach).

 \textbf{Postulate 1.}
\textit{Existence of infinite and infinitesimal objects is
postulated but it is also accepted  that human beings and machines
are able to execute only a finite number of operations.}

\textbf{Postulate 2.} \textit{It is not discussed \textbf{what
are} the mathematical objects we deal with; we just construct more
powerful tools that   allow us to improve our capacities to
observe and to describe properties of mathematical objects.}

\textbf{Postulate 3.} \textit{The principle formulated by Ancient
Greeks  `The part is less than the whole' is applied to all
numbers (finite, infinite, and infinitesimal) and to all sets and
processes (finite and infinite).}

Due to this declared applied statement, such traditional concepts
as bijection, numerable and continuum sets, cardinal and ordinal
numbers are not applied when one works with the Infinity Computer
because they belong to Cantor's approach having significantly more
theoretical character and based on   different assumptions.
However, the methodology used by the Infinity Computer does not
contradict Cantor. In contrast, it evolves his deep ideas
regarding existence of different infinite numbers in a more
practical way.

The accepted applied methodology means in particular  (see
Postulate~1) that we shall never be able to give a complete
description of infinite processes and sets due to our finite
capabilities. Acceptance of Postulate~1 means also that we
understand that we are able to write down only a finite number of
symbols to express numbers in any numeral system\footnote{We
remind that a \textit{numeral}  is a symbol or group of symbols
that represents a \textit{number}. A \textit{number} is a concept
that a \textit{numeral} expresses. The same number can be
represented by different numerals (e.g., the symbols `10', `ten',
and `X' are different numerals, but they all represent the same
number).}.

Postulate 2 states that the philosophical triad -- researcher,
object of investigation, and tools used to observe the object --
existing in such natural sciences as Physics and Chemistry, exists
in Mathematics, too. In natural sciences, the instrument used to
observe the object limits and influences the results of
observations. The same happens in Mathematics where numeral
systems used to express numbers are among the instruments of
observations used by mathematicians. Usage of a powerful numeral
system  gives a possibility to obtain more precise results in
Mathematics in the same way as usage of a good microscope gives a
possibility to obtain more precise observations in Physics.
However, due to Postulate~1, the capabilities of the tools will be
always limited.

Particularly, this means that from this applied point of view,
axiomatic systems do not define mathematical objects but just
determine formal rules for operating with certain numerals
reflecting some properties of the studied mathematical objects.
For example,  axioms for real numbers are considered together with
a particular numeral system $\mathcal{S}$ used to write down
numerals and are viewed as practical rules (associative and
commutative properties of multiplication and addition,
distributive property of multiplication over addition, etc.)
describing possible operations with the numerals. The completeness
property is interpreted as a possibility to extend $\mathcal{S}$
with additional symbols (e.g., $e$, $\pi$, $\sqrt{2}$, etc.)
taking care of the fact that the results of computations with
these symbols  agree with facts observed in practice. As a rule,
assertions regarding numbers that cannot be expressed in a numeral
system are avoided (e.g., it is not supposed that real numbers
form a field).

In this paper,  the methodology from
\cite{Sergeyev,Poland,chaos,spirals,informatica} is used to
describe how the Infinity Computer can be applied for solving new
and old (but with higher precision) computational problems.
Representation of infinite and infinitesimal numbers  at the
Infinity Computer  and operations with them are described in
Section~\ref{s1}. Then, Sections~\ref{s4} and~\ref{s_calcolo4}
present results dealing with applications related to linear
algebra and calculating divergent series, limits, volumes, and
probabilities.

\section{Representation of infinite and  infinitesimal numbers at
the Infinity Computer} \label{s1}

In \cite{Sergeyev,Poland,chaos,spirals,informatica},   a new
powerful numeral system has been developed to express finite,
infinite, and infinitesimal numbers in a unique framework by a
finite number of symbols. The main idea of the new approach
consists of measuring infinite and infinitesimal quantities   by
different (infinite, finite, and infinitesimal) units of measure.
This section gives just a brief tour to the representation of
infinite and infinitesimal numbers at the Infinity Computer and
describes how operations with them can be executed. It allows us
to introduce the necessary notions and designations. In order to
have a comprehensive introduction to the new methodology, we
invite the reader to have a look at the recent survey
\cite{informatica} downloadable from \cite{www} or at the book
\cite{Sergeyev} (written in a popular manner) before approaching
Sections~\ref{s4} and~\ref{s_calcolo4}.

\subsection{A new infinite numeral and a positional numeral system with\\ the~infinite radix}

A new infinite unit of measure   has been introduced for this
purpose   as the number of elements of the set $\mathbb{N}$ of
natural numbers. It is expressed by a new numeral \ding{172}
called \textit{grossone}. It is necessary to emphasize immediately
that the infinite number \ding{172} is not either Cantor's
$\aleph_0$ or $\omega$. In particular, \ding{172} has both
cardinal and ordinal properties as usual finite natural numbers.
Formally, grossone is introduced as a new number by describing its
properties postulated by the \textit{Infinite Unit Axiom} (IUA)
(see \cite{Sergeyev,www,Mathesis,chaos,informatica}). This axiom
is added to axioms for real numbers (viewed in sense of
Postulates~1--3) similarly to addition of the axiom determining
zero to axioms of natural numbers when integer numbers are
introduced.

One of the important differences of the new approach with respect
to the non-standard analysis consists of the fact that
$\mbox{\ding{172}} \in \mathbb{N}$ because grossone has been
introduced as the quantity of natural numbers  (similarly, the
number 3 being the number of elements of the set $\{1, 2, 3 \}$ is
the largest element in this set). The new numeral \ding{172}
allows one to write down the set, $\mathbb{N}$, of natural numbers
in the form
 \beq
\mathbb{N} = \{ 1,2,3, \hspace{5mm} \ldots  \hspace{5mm}
 \mbox{\ding{172}}-3, \hspace{2mm} \mbox{\ding{172}}-2,
\hspace{2mm}\mbox{\ding{172}}-1, \hspace{2mm} \mbox{\ding{172}} \}
\label{4.1}
       \eeq
   where the numerals
 \beq
\ldots  \hspace{2mm} \mbox{\ding{172}}-3,
\hspace{2mm}\mbox{\ding{172}}-2, \hspace{2mm}\mbox{\ding{172}}-1,
\hspace{2mm} \mbox{\ding{172}} \label{4.2}
       \eeq
 indicate \textit{infinite} natural numbers.

It is important to emphasize that in the new approach the set
(\ref{4.1}) is the same set of natural numbers
 \beq
\mathbb{N} = \{ 1,2,3, \hspace{2mm} \ldots \hspace{2mm}  \}
\label{4.1_calcolo}
       \eeq
 we
are used to deal with and infinite numbers (\ref{4.2}) also belong
to  $\mathbb{N}$. Both records, (\ref{4.1}) and
(\ref{4.1_calcolo}), are correct and do not contradict each other.
They just use two different numeral systems to express
$\mathbb{N}$. Traditional numeral systems   do not allow us to see
infinite natural numbers that we can observe now thanks
to~\ding{172}. Similarly, a primitive tribe of Pirah\~{a} (see
\cite{Gordon}) living in Amazonia   that uses a very weak numeral
system for counting (one, two, many)  is not able to see finite
natural numbers greater than~2. In spite of this fact, these
numbers (e.g., 3 and 4) belong to $\mathbb{N}$ and are visible if
one uses a more powerful numeral system. Thus, we have the same
object of observation -- the set $\mathbb{N}$ -- that can be
observed by different instruments -- numeral systems -- with
different accuracies (see Postulate~2).

It is worthy to notice that the introduction of \ding{172}
defines the set of \textit{extended natural numbers} indicated as
$\widehat{\mathbb{N}}$ and including $\mathbb{N}$ as a proper
subset
 \beq
  \widehat{\mathbb{N}} = \{
1,2, \ldots ,\mbox{\ding{172}}-1, \mbox{\ding{172}},
\mbox{\ding{172}}+1, \ldots , \mbox{\ding{172}}^2-1,
\mbox{\ding{172}}^2, \mbox{\ding{172}}^2+1, \ldots \}.
\label{4.2.2}
       \eeq
Due to Postulates 1 and 2, the new numeral system cannot give
answers to \textit{all} questions regarding infinite sets. What
can we say, for instance, about the number of elements of the set
$\widehat{\mathbb{N}}$? The introduced numeral system based on
\ding{172} is too weak to give an answer  to this question. It is
necessary to introduce in a way a more powerful numeral system by
defining new numerals (for instance, \ding{173}, \ding{174}, etc).

Inasmuch as it has been postulated that grossone is a number,
associative and commutative properties of multiplication and
addition, distributive property of multiplication over addition,
existence of   inverse  elements with respect to addition and
multiplication hold for grossone as for finite numbers.
Particularly, this means that  the following relations hold for
grossone, as for any other number
 \beq
 0 \cdot \mbox{\ding{172}} =
\mbox{\ding{172}} \cdot 0 = 0, \hspace{3mm}
\mbox{\ding{172}}-\mbox{\ding{172}}= 0,\hspace{3mm}
\frac{\mbox{\ding{172}}}{\mbox{\ding{172}}}=1, \hspace{3mm}
\mbox{\ding{172}}^0=1, \hspace{3mm}
1^{\mbox{\tiny{\ding{172}}}}=1, \hspace{3mm}
0^{\mbox{\tiny{\ding{172}}}}=0.
 \label{3.2.1}
       \eeq

To express infinite and infinitesimal numbers  at the Infinity
Computer, records similar to traditional positional numeral
systems can be used (see \cite{Sergeyev,www,Mathesis,chaos}). In
order to construct a number $C$ in the new numeral positional
system with the radix \ding{172}, we subdivide $C$ into groups
corresponding to powers of \ding{172}:
 \beq
  C = c_{p_{m}}
\mbox{\ding{172}}^{p_{m}} +  \ldots + c_{p_{1}}
\mbox{\ding{172}}^{p_{1}} +c_{p_{0}} \mbox{\ding{172}}^{p_{0}} +
c_{p_{-1}} \mbox{\ding{172}}^{p_{-1}}   + \ldots   + c_{p_{-k}}
 \mbox{\ding{172}}^{p_{-k}}.
\label{3.12}
       \eeq
 Then, the record
 \beq
  C = c_{p_{m}}
\mbox{\ding{172}}^{p_{m}}    \ldots   c_{p_{1}}
\mbox{\ding{172}}^{p_{1}} c_{p_{0}} \mbox{\ding{172}}^{p_{0}}
c_{p_{-1}} \mbox{\ding{172}}^{p_{-1}}     \ldots c_{p_{-k}}
 \mbox{\ding{172}}^{p_{-k}}
 \label{3.13}
       \eeq
represents  the number $C$, where  $c_i$ are called
\textit{grossdigits} and are expressed by traditional numeral
systems used to represent finite numbers (e.g., floating point
numerals). Grossdigits  can be both positive and negative. They
show how many corresponding units should be added or subtracted in
order to form the number $C$. Grossdigits can be expressed by
several symbols.

Numbers $p_i$ in (\ref{3.13}) called \textit{grosspowers}  can be
finite, infinite, and infinitesimal. They   are sorted in
decreasing order with $ p_0=0$:
\[
p_{m} >  p_{m-1}  > \ldots    > p_{1} > p_0 > p_{-1}  > \ldots
p_{-(k-1)}  >   p_{-k}.
 \]
 In the record (\ref{3.13}), we write
$\mbox{\ding{172}}^{p_{i}}$ explicitly because in the new numeral
positional system  the number   $i$ in general is not equal to the
grosspower $p_{i}$.

\textit{Finite numbers} in this new numeral system are represented
by numerals having only the grosspower $ p_0=0$. In fact, if we
have a number $C$ such that $m=k=$~0 in representation
(\ref{3.13}), then due to (\ref{3.2.1}),   we have $C=c_0
\mbox{\ding{172}}^0=c_0$. Thus, the number $C$ in this case does
not contain grossone and is equal to the grossdigit $c_0$ being a
conventional finite number     expressed in a traditional finite
numeral system.

\textit{Infinitesimal numbers}   are represented by numerals $C$
having only negative finite or infinite grosspowers, e.g.,
$3.48\mbox{\ding{172}}^{-46.71}26.4\mbox{\ding{172}}^{-132\mbox{\tiny{\ding{172}}}}$.
The simplest infinitesimal number is
$\frac{1}{\mbox{\ding{172}}}=\mbox{\ding{172}}^{-1}$ being the
inverse element with respect to multiplication for \ding{172}:
 \beq
\mbox{\ding{172}}^{-1}\cdot\mbox{\ding{172}}=\mbox{\ding{172}}\cdot\mbox{\ding{172}}^{-1}=1.
 \label{3.15.1}
       \eeq
Note that all infinitesimals are not equal to zero. Particularly,
$\frac{1}{\mbox{\ding{172}}}>0$ because it is a result of division
of two positive numbers.

\textit{Infinite numbers}    are expressed by numerals having at
least one positive finite or infinite grosspower. Thus, they have
infinite parts and can also have a finite part and infinitesimal
ones.  For instance, the number
\[
 12.4\mbox{\ding{172}}^{34.21\mbox{\tiny{\ding{172}}}}(\mbox{\small-}20.64)\mbox{\ding{172}}^{15}
0.8\mbox{\ding{172}}^{0}0.71\mbox{\ding{172}}^{-3}32.1\mbox{\ding{172}}^{-6.5\mbox{\tiny{\ding{172}}}}
\]
has two infinite parts,
$12.4\mbox{\ding{172}}^{34.21\mbox{\tiny{\ding{172}}}}$ and
$-20.64\mbox{\ding{172}}^{15}$, one finite part,
$0.8\mbox{\ding{172}}^{0}=0.8$, and two infinitesimal parts,
$0.71\mbox{\ding{172}}^{-3}$ and
$32.1\mbox{\ding{172}}^{-6.5\mbox{\tiny{\ding{172}}}}$.

\subsection{Arithmetical operations executed by the Infinity Computer}

A working software simulator of the Infinity Computer has been
implemented (see \cite{Sergeyev_patent,www,Poland}). It works with
infinite, finite, and infinitesimal numbers \textit{numerically},
(not symbolically) and executes the arithmetical operations as
follows.

Let us consider the operation of \textit{addition}
(\textit{subtraction} is a direct consequence of addition and is
thus omitted) of two given infinite numbers $A$ and $B$, where
 \beq
A= \sum_{i=1}^{K} a_{k_{i}}\mbox{\ding{172}}^{k_{i}}, \hspace{1cm}
B= \sum_{j=1}^{M} b_{m_{j}}\mbox{\ding{172}}^{m_{j}}, \hspace{1cm}
C= \sum_{i=1}^{L} c_{l_{i}}\mbox{\ding{172}}^{l_{i}},
 \label{3.20}
       \eeq
and the result $C=A+B$ is constructed    by including in it all
items $a_{k_{i}}\mbox{\ding{172}}^{k_{i}}$ from $A$ such that
$k_{i} \neq m_{j},1 \le j \le M,$ and all items
$b_{m_{j}}\mbox{\ding{172}}^{m_{j}}$ from $B$ such that $m_{j}
\neq k_{i},1 \le i \le K$. If in $A$ and $B$ there are items such
that $k_{i}=m_{j}$, for some $i$ and $j$, then this grosspower
$k_{i}$ is included in $C$ with the grossdigit
$b_{k_{i}}+a_{k_{i}}$, i.e., as
$(b_{k_{i}}+a_{k_{i}})\mbox{\ding{172}}^{k_{i}}$.

 The operation of \textit{multiplication}  of two numbers
 $A$ and $B$ in the form (\ref{3.20}) returns, as
the result, the infinite number $C$ constructed as follows:
 \beq
C= \sum_{j=1}^{M} C_{j}, \hspace{5mm} C_{j} =
b_{m_{j}}\mbox{\ding{172}}^{m_{j}}\cdot A =\sum_{i=1}^{K}
a_{k_{i}}b_{m_{j}}\mbox{\ding{172}}^{k_{i}+m_{j}}, \hspace{5mm}1
\le j \le M. \label{3.23}
       \eeq

In the operation of \textit{division} of a   number $C$ by a
number $B$ from (\ref{3.20}), we obtain a result $A$ and a
reminder $R$ (that can be also equal to zero), i.e., $C =A \cdot
B+R$. The number $A$ is constructed as follows.  The first
grossdigit $a_{k_{K}}$ and the corresponding maximal exponent
${k_{K}}$ are established from the equalities
 \beq
a_{k_{K}}=c_{l_{L}}/ b_{m_{M}}, \hspace{4mm}  k_{K} = l_{L}-
m_{M}.
 \label{3.25}
       \eeq
Then the first partial reminder  $R_1$ is calculated as
 \beq
R_1= C - a_{k_{K}}\mbox{\ding{172}}^{k_{K}} \cdot B.
\label{3.25.0}
       \eeq
If $R_1 \neq 0$ then the number $C$ is substituted by $R_1$ and
the process is repeated by a complete analogy. The process stops
when a partial reminder equal to zero is found (this means that
the final reminder $R=0$) or when a required accuracy of the
result is reached.

\begin{example}
\label{e6} We consider two infinite numbers $A$ and $B$, where
$$
A=304.21\mbox{\ding{172}}^{16.8\mbox{\tiny{\ding{172}}}}(\mbox{\small-}7.1)\mbox{\ding{172}}^{12}
41.2\mbox{\ding{172}}^{0}, \hspace{1cm}
B=6.23\mbox{\ding{172}}^{3}
13.1\mbox{\ding{172}}^{0}15\mbox{\ding{172}}^{-6.2\mbox{\tiny{\ding{172}}}}.
$$
Their sum $C$ is calculated as follows:
\[
C=A+B=304.21\mbox{\ding{172}}^{16.8\mbox{\tiny{\ding{172}}}}+(\mbox{\small-}7.1)\mbox{\ding{172}}^{12}+
41.2\mbox{\ding{172}}^{0}+ 6.23\mbox{\ding{172}}^{3}+
13.1\mbox{\ding{172}}^{0}
+15\mbox{\ding{172}}^{-6.2\mbox{\tiny{\ding{172}}}}=
\]
\[
304.21\mbox{\ding{172}}^{16.8\mbox{\tiny{\ding{172}}}}-
7.1\mbox{\ding{172}}^{12}+6.23\mbox{\ding{172}}^{3}+
54.3\mbox{\ding{172}}^{0}+15\mbox{\ding{172}}^{-6.2\mbox{\tiny{\ding{172}}}}=
 \]
\[
\begin{tabular}{cr}\hspace {28mm}$304.21\mbox{\ding{172}}^{16.8\mbox{\tiny{\ding{172}}}}
(\mbox{\small-}7.1)\mbox{\ding{172}}^{12}6.23\mbox{\ding{172}}^{3}
54.3\mbox{\ding{172}}^{0}15\mbox{\ding{172}}^{-6.2\mbox{\tiny{\ding{172}}}}.\hfill
{ }$ & \hspace {14mm}
 $\Box$
\end{tabular}
\]
\end{example}
More examples illustrating the work of the Infinity Computer  can
be found in \cite{Sergeyev,Poland,Mathesis,chaos,informatica}.

\section{Examples of situations where the Infinity Computer\\ executes
operations that traditionally   required  a human intervention}
\label{s4}

In this section, we describe a number of   computational tools
provided by the new methodology and the Infinity Computer. It
becomes possible in several occasions to automatize the process of
the solving of computational problems avoiding an interruption of
the work of computer procedures and the necessity of
 a human intervention required when one works with traditional computers. For
instance, when one meets a necessity to work with divergent series
or, even worse, their difference, traditional computers are not
able to execute these operations automatically and a human
intervention is required in order to avoid these difficulties. In
this situation and in other examples below, it is shown how the
work usually done by humans can be formalized following
Postulates~1--3 and passed to the Infinity Computer.

It is necessary to emphasize that the examples described in this
section are related to \textit{numerical} computations at the
Infinity Computer. No symbolic computations  are required to work
with infinite and infinitesimal numbers when one uses the Infinity
Computer.

\subsection{Calculating sums with an infinite number of items}
\label{s4.3}

The new approach allows one to use the Infinity Computer for the
calculation of sums with an infinite number of items. The term
`series' is not used here because, due to Postulate~3,  it is
required to indicate explicitly the number of items (finite or
infinite) in any sum. Naturally, it is necessary that the number
of items and the result of the considered sum are expressible in
the numeral system used for calculations.

Let us illustrate the new possibilities by considering   two
traditional infinite series $S_1=1+1+1+\ldots$ and
$S_2=30+30+30+\ldots$  The traditional analysis gives us a very
poor answer that both of them diverge to infinity and, therefore,
the results cannot be calculated and represented at the
traditional computers. Such operations as $S_2 - S_1$ or $
\frac{S_1}{S_2} $ are not defined and humans should return to the
original  physical problem in order to understand whether there
exist answers to these questions.

Now, when we are able to express not only different finite numbers
but also different infinite numbers, it is necessary to indicate
explicitly the number of items in the sums $S_1$ and $S_2$ and it
is not important whether it is finite or infinite.  Due to
Postulate~3, by changing the number of items in the sums, we
change the respective results, too.

Suppose that the   sum $S_1$ has $k$ items and the sum $S_2$ has
$n$ items. Then
$$S_1(k)=\underbrace{1+1+1+\ldots+1}_k, \hspace{1cm} S_2(n)=\underbrace{30+30+30+\ldots+30}_n.$$
and it follows $S_1(k)=k$ and $S_2(n)=30n$. If the   results
$S_1(k)=k$ and $S_2(n)=30n$ are expressible in the chosen numeral
system, then indeterminate  forms disappear and the expressions
$S_2(k) - S_1(n)$ and $ \frac{S_1(k)}{S_2(n)}$ can be easily
calculated using the Infinity Computer.

If, for instance, $k=n=5\mbox{\ding{172}}$ then we obtain
$S_1(5\mbox{\ding{172}})=5\mbox{\ding{172}}$,
$S_2(5\mbox{\ding{172}})=150\mbox{\ding{172}}$ and
 $S_2(5\mbox{\ding{172}}) -  S_1(5\mbox{\ding{172}})
= 145\mbox{\ding{172}} > 0$.

If  $k=30\mbox{\ding{172}}$ and $n=\mbox{\ding{172}}$   we obtain
$S_1(30\mbox{\ding{172}})=30\mbox{\ding{172}}$,
$S_2(\mbox{\ding{172}})=30\mbox{\ding{172}}$ and  it follows
 $S_2(\mbox{\ding{172}}) -
S_1(30\mbox{\ding{172}})=0$.

If   $k=30\mbox{\ding{172}}2$ (we remind that we use here a
shorten way to write down this infinite number, the complete
record is $30\mbox{\ding{172}}^{1}2\mbox{\ding{172}}^{0}$) and
$n=\mbox{\ding{172}}$ we obtain
$S_1(30\mbox{\ding{172}}2)=30\mbox{\ding{172}}2$,
$S_2(\mbox{\ding{172}})=30\mbox{\ding{172}}$ and  it follows
 \[
 S_2(\mbox{\ding{172}}) -
S_1(30\mbox{\ding{172}}2)= 30\mbox{\ding{172}} -
30\mbox{\ding{172}}2= - 2< 0.
 \]
  Analogously, the expression $
\frac{S_1(k)}{S_2(n)} $ can be calculated.

Let us consider now  the famous divergent series with alternate
signs $S_3=1-1+1-1+\ldots$ In literature there exist  many
approaches giving different answers regarding the value of this
series (see \cite{Knopp}). All of them use various notions of
average. However, the notions of sum and average are different. In
our approach, we do not appeal to average and calculate the
required sum directly. To do this we should indicate explicitly
the number of items, $k$, in the sum. Then
\[
S_3(k)=\underbrace{1-1+1-1+1-1+1-\ldots}_{k} = \left \{
\begin{array}{ll} 0, &
  \mbox{if  } k=2n,\\
1, &    \mbox{if  } k=2n+1,\\
 \end{array} \right.
\]
and it is not important whether $k$ is  finite or infinite.

We conclude this subsection by studying the
 series $S_4 =\sum_{i=1}^{\infty}\frac{1}{2^i}$   converging to one.
 The new approach allows us  to give a more precise answer. Due to Postulate~3, the
formula
$$S_4(k)=\sum_{i=1}^{k}\frac{1}{2^i}=1-\frac{1}{2^k}$$
can be used directly for   infinite $k$, too. For example, if
$k=\mbox{3\ding{172}}$ then
$$S_4(3\mbox{\ding{172}})=\sum_{i=1}^{3\mbox{\tiny{\ding{172}}}}\frac{1}{2^i}=1-\frac{1}{2^{3\mbox{\tiny{\ding{172}}}}},$$
where $\frac{1}{2^{3\mbox{\tiny{\ding{172}}}}}$ is infinitesimal.
Thus, the traditional answer $\sum_{i=1}^{\infty}\frac{1}{2^i}=1$
is just a finite approximation to our more precise result using
infinitesimals. The traditional numeral systems do not allow us to
distinguish results of the sums for infinite values of $k$.  More
examples can be found in \cite{chaos}.

Thus, if one is able to calculate a partial sum of a series $S$,
he/she can use the  formula applied for this calculation to
evaluate at the Infinity Computer sums $S(k)$    with $k$ items
for finite and infinite values of $k$ and finite, infinite, and
infinitesimal values of $S(k)$ and to use the obtained results in
further calculations.

\subsection{Computing expressions with infinite and infinitesimal arguments}
\label{s4.5}

In the traditional analysis, the concept of the limit has been
introduced in order to avoid difficulties that one faces when
he/she wants to evaluate an expression at infinity or at a point
$x$ infinitely close to a point $a$. If $\lim_{x \rightarrow
a}f(x)$ exists, then it gives us a very poor -- just one  value --
information about the behavior of $f(x)$ when $x$ tends to $a$.

Now we can obtain significantly more rich information using the
Infinity Computer independently on the fact of existence of the
limit. We can calculate $f(x)$ directly at any finite, infinite,
or infinitesimal point expressible in the new positional system
even if the limit does not exist. Thus, limits   can be
substituted by precise numerals $f(a)$ that are different for
different infinite, finite, or infinitesimal values of $x=a$. This
is very important for practical computations because this
substitution eliminates indeterminate forms, i.e., again the
Infinity Computer should not stop its calculations as traditional
computers are forced to do when they encounter indeterminate
forms.

\begin{example}
\label{e15} In the traditional analysis, the following two limits
 \[
  \lim_{x \rightarrow +\infty}(x^4+11.5x^2+10^{100})= +\infty, \hspace{1cm}
 \lim_{x \rightarrow +\infty}(x^4+11.5x^2)= +\infty.
 \]
give  us   the same result, $+\infty$,  in spite of the fact that
for any finite $x$ the  difference between the two expressions is
equal to quite a large number
 \[
x^4+11.5x^2+10^{100} - (x^4+11.5x^2) =  10^{100}.
\]
The new approach allows us to calculate exact values of both
expressions, $x^4+11.5x^2+10^{100}$ and $x^4+11.5x^2$, at any
infinite (and infinitesimal) $x$ expressible in the chosen numeral
system. For instance, the choice of $x=\mbox{3\ding{172}}^{2}$
gives the values
\[
(3\mbox{\ding{172}}^{2})^{4}+11.5(3\mbox{\ding{172}}^{2})^{2}+10^{100}=
81\mbox{\ding{172}}^{8}103.5\mbox{\ding{172}}^{4}10^{100}\mbox{\ding{172}}^{0}
\]
and $81\mbox{\ding{172}}^{8}103.5\mbox{\ding{172}}^{4}$,
respectively. Consequently, one obtains
 \[
\begin{tabular}{cr}\hspace {24mm}$
81\mbox{\ding{172}}^{8}103.5\mbox{\ding{172}}^{4}10^{100}\mbox{\ding{172}}^{0}
- 81\mbox{\ding{172}}^{8}103.5\mbox{\ding{172}}^{4} = 10^{100}.$ &
\hspace {17mm}
 $\Box$
\end{tabular}
\]
\end{example}

An additional advantage of the usage of the Infinity Computer
arises in the following situations. Suppose that   we have a
computer procedure calculating $f(x)$, we do not know the
corresponding analytic formulae for $f(x)$, for a certain argument
$a$ the value $f(a)$ is not defined (or a traditional computer
produces an overflow or underflow message), and it is necessary to
calculate the $\lim_{x \rightarrow a}f(x)$. Traditionally, this
situation requires a human intervention and an additional
theoretical investigation whereas the Infinity Computer is able to
process it automatically working numerically with the expressions
involved in the procedure. It is sufficient to calculate $f(x)$,
for example, at a point $x=a+\mbox{\ding{172}}^{-1}$ in cases of
finite $a$ or $a=0$ and $x=\mbox{\ding{172}}$ in the case when we
are interested in the behavior of $f(x)$ at infinity.

\begin{example}
\label{e_calcolo_1} Suppose that we have two procedures evaluating
$f(x)=\frac{x^2+2x}{x}$ and  $g(x)=\frac{34}{x}$. Obviously,
$f(0)$ and $g(0)$ are not defined and it is not possible to
calculate $\lim_{x \rightarrow 0}f(x)$, $\lim_{x \rightarrow
\infty}f(x)$ and $\lim_{x \rightarrow 0}g(x)$, $\lim_{x
\rightarrow \infty}g(x)$ using traditional computers. Then,
suppose that we are interested in evaluating the expression
\[
h(x)= (f(x)-2)\cdot g(x).
\]
It is easy to see that $h(x)=34$ for any finite value of $x$. On
the other hand, the following limits
\[
\lim_{x \rightarrow 0}h (x)= (\lim_{x \rightarrow 0}f(x)-2)\cdot
\lim_{x \rightarrow 0}g(x),
\]
\[
\lim_{x \rightarrow \infty}h (x)= (\lim_{x \rightarrow
\infty}f(x)-2)\cdot \lim_{x \rightarrow \infty}g(x)
\]
cannot be evaluated. The Infinity Computer can calculate $h(x)$
numerically for different  infinitesimal and infinite values of
$x$ obtaining the same result that takes place for finite $x$. For
example, it follows
\[
h (\mbox{\ding{172}}^{-1})=
\left(\frac{(\mbox{\ding{172}}^{-1})^2+2\mbox{\ding{172}}^{-1}}{\mbox{\ding{172}}^{-1}}-2\right)\cdot
\frac{34}{\mbox{\ding{172}}^{-1}}=(\mbox{\ding{172}}^{-1}+2-2)\cdot
34\mbox{\ding{172}}=34,
\]
\[
\hspace{20mm}h (\mbox{\ding{172}})=
\left(\frac{\mbox{\ding{172}}^2+2\mbox{\ding{172}}}{\mbox{\ding{172}}}-2\right)\cdot
\frac{34}{\mbox{\ding{172}}}=(\mbox{\ding{172}}+2-2)\cdot
34\mbox{\ding{172}}^{-1}=34. \hspace{15mm}  \Box
\]
\end{example}

It is worthy to notice  that expressions can be calculated by the
Infinity Computer for infinite and infinitesimal arguments, even
when their limits do not exist, thus giving     a very powerful
tool for studying  divergent processes.
\begin{example}
\label{e16}
 The   limit
$ \lim_{n \rightarrow +\infty}f(n),$ where $f(n)=(-5)^n n,$ does
not exist. However, we can calculate the expression $(-5)^n n$ for
different infinite values of $n$. For instance, since it can be
easily proved that grossone is even (see
\cite{Sergeyev,informatica}), for $n=\mbox{\ding{172}}$ it follows
$f(\mbox{\ding{172}})=5^{{\mbox{\tiny\ding{172}}}}\mbox{\ding{172}}$
and for $n=\mbox{\ding{172}}-1$ we have
$f(\mbox{\ding{172}}-1)=-5^{{\mbox{\tiny\ding{172}}}-1}(\mbox{\ding{172}}-1)$.
\hfill $\Box$
\end{example}

Thus, these new computational possibilities of the Infinity
Computer allow one both to avoid calculating limits theoretically
and to increase the accuracy of numerical computations.

\subsection{ Usage of infinitesimals  for solving systems of linear equations}
\label{s_calcolo3.3}

Very often in computations,   an algorithm  performing
calculations
 encounters a situation where the problem to divide by zero
occurs.  Then, obviously,  this operation cannot be executed. If
it is known that the problem under consideration has a solution,
then
 a number of additional computational steps trying to avoid  this
division is performed. A typical example of this kind is the
operation of pivoting used when one solves systems of linear
equations by  an algorithm such as Gauss-Jordan elimination.
Pivoting is the interchanging of rows (or both rows and columns)
in order to avoid division by zero and to place a particularly
`good' element in the diagonal position prior to a particular
operation.

The following two simple examples give just an idea of a numerical
usage of infinitesimals and show that the usage of infinitesimals
can help to avoid pivoting in cases when the pivotal element is
equal to zero. We emphasize again that the Infinity Computer (see
\cite{Sergeyev_patent}) works with infinite and infinitesimal
numbers expressed in the positional numeral system (\ref{3.12}),
(\ref{3.13}) numerically, not symbolically.

\begin{example}
\label{e_calcolo_2}

Solution to the system
\[
 \begin{array}{ccc}
 \left[ \begin{array}{cc}  0 & 1 \\ 2 & 2 \end{array}  \right]  &
 \left[ \begin{array}{c}  x_1 \\ x_2 \end{array} \right] = &
  \left[ \begin{array}{c}  2 \\ 2 \end{array} \right]
  \end{array}
\]
is obviously given by   $x^*_1=  -1,$ $ x^*_2=2$. It cannot be
found by the method of Gauss without pivoting since   the first
pivotal element $a_{11}= 0$.

Since all the elements of the matrix are finite numbers, let us
substitute the element $a_{11}=0$ by $\mbox{\ding{172}}^{-1}$ and
perform exact  Gauss transformations without pivoting:
\[
 \left[ \begin{array}{cc|c}   \mbox{\ding{172}}^{-1}  & 1 & 2 \\ 2 & 2 & 2\end{array}  \right]
 \rightarrow
 \left[ \begin{array}{cc|c}  1 & \mbox{\ding{172}} & 2\mbox{\ding{172}} \\ 0 & -2\mbox{\ding{172}}+2 & -4\mbox{\ding{172}}+2 \end{array}  \right]
 \rightarrow
 \left[ \begin{array}{cc|c}  1 & \mbox{\ding{172}} & 2\mbox{\ding{172}} \\ 0 & 1 &  \frac{-4\mbox{\small{\ding{172}}}+2}{-2\mbox{\small{\ding{172}}}+2}\end{array}  \right]
\]
\[
 \left[ \begin{array}{cc|c}  1 & 0 & 2\mbox{\ding{172}}- \mbox{\ding{172}}\cdot\frac{-4\mbox{\small{\ding{172}}}+2}{-2\mbox{\small{\ding{172}}}+2}  \\ 0 & 1 & \frac{-4\mbox{\small{\ding{172}}}+2}{-2\mbox{\small{\ding{172}}}+2} \end{array}
 \right]
 \rightarrow
 \left[ \begin{array}{cc|c}  1 & 0 &  \frac{2\mbox{\small{\ding{172}}}}{-2\mbox{\small{\ding{172}}}+2}  \vspace*{1mm}\\ 0 & 1 & \frac{-4\mbox{\small{\ding{172}}}+2}{-2\mbox{\small{\ding{172}}}+2} \end{array}
 \right]
  \rightarrow
 \left[ \begin{array}{cc|c}  1 & 0 &   -1 + \frac{1}{1-\mbox{\small{\ding{172}}}}  \\ 0 & 1 &  2-\frac{1}{1-\mbox{\small{\ding{172}}}}  \end{array}
 \right].
\]
It follows immediately that the solution to the initial system is
given by the finite parts  of  numbers $ -1 +
\frac{1}{1-\mbox{\small{\ding{172}}}} $ and $
2-\frac{1}{1-\mbox{\small{\ding{172}}}}$.

We have introduced the number $\mbox{\ding{172}}^{-1}$ once and,
as a result, we have obtained expressions where the maximal power
of grossone is one and there are rational expressions depending on
grossone, as well. It is possible to manage these rational
expressions in two ways: (i) to execute division in order to
obtain   its result in the form (\ref{3.12}), (\ref{3.13}); (ii)
without executing division. In the latter case, we just continue
to work with rational expressions. In the case (i), since  we need
finite numbers as final results,   in the  result of   division it
is not necessary to store the parts $c_p\mbox{\ding{172}}^p$ with
$p < -1$. These parts can be forgotten because in any way the
result of their successive multiplication with the numbers of the
type $c_1\mbox{\ding{172}}^1$ (remind that 1 is the maximal
exponent present in the matrix under consideration) will give
exponents less than zero, i.e., numbers with these exponents will
be infinitesimals that are not interesting for us in this
computational context.

Thus, by using the positional numeral system (\ref{3.12}),
(\ref{3.13}) with the radix grossone we obtain
\[
 \left[ \begin{array}{cc|c}  1 & \mbox{\ding{172}} & 2\mbox{\ding{172}} \\ 0 & 1 &  \frac{-4\mbox{\small{\ding{172}}}+2}{-2\mbox{\small{\ding{172}}}+2}\end{array}  \right]
  \rightarrow
\left[ \begin{array}{cc|c}  1 & \mbox{\ding{172}} &
2\mbox{\ding{172}} \\ 0 & 1 &
 2\mbox{\ding{172}}^0\mbox{\small{+}}1\mbox{\ding{172}}^{-1}\end{array}
\right]
  \]
\[
\left[ \begin{array}{cc|c}  1 & 0 &
2\mbox{\ding{172}}-\mbox{\ding{172}}\cdot ( 2\mbox{\ding{172}}^0
1\mbox{\ding{172}}^{-1})
\\ 0 & 1 &
 2\mbox{\ding{172}}^01\mbox{\ding{172}}^{-1}\end{array}
\right]
  \rightarrow
\left[ \begin{array}{cc|c}  1 & 0 & -1\mbox{\ding{172}}^0\\
 0 & 1 &
 2\mbox{\ding{172}}^0 1\mbox{\ding{172}}^{-1}\end{array}
\right].
\]
The finite parts  of  numbers $-1\mbox{\ding{172}}^0$ and
$2\mbox{\ding{172}}^0 1\mbox{\ding{172}}^{-1}$, i.e., $-1$ and 2
respectively, then provide  the required solution. \hfill $\Box$
\end{example}

%The next example shows that if the second pivotal element $a_{jj}=
%0$ is encountered, then  to avoid pivoting  it is necessary   to
%substitute it by $\mbox{\ding{172}}^{-2}$.

\begin{example}
\label{e_calcolo_3}

Solution to the system
\[
 \begin{array}{ccc}
 \left[ \begin{array}{ccr}  0 & 0 & 1  \\ 2 & 0 & -1 \\ 1 & 2 & 3
\end{array}  \right]  &
 \left[ \begin{array}{c}  x_1 \\ x_2 \\ x_3 \end{array} \right] = &
  \left[ \begin{array}{c}  1 \\  3 \\ 1 \end{array} \right]
  \end{array}
\]
is the following: $ x^*_1= 2, $ $    x^*_2  = -2, $ and $   x^*_3
= 1$. The coefficient matrix of this system has the first two
leading principal minors equal to zero. Consequently, the first
two pivots, in the Gauss transformations, are zero. We solve the
system without pivoting by substituting  the zero pivot by
$\mbox{\ding{172}}^{-1}$, when necessary.

Let us show how the exact computations are executed:
\[
  \left[ \begin{array}{ccr|c}  0 & 0  & 1 & 1 \\ 2 & 0 & -1
& 3 \\ 1 & 2 & 3 & 1 \end{array}  \right]
 \rightarrow
 \left[ \begin{array}{ccc|c}  1 & 0  & \mbox{\ding{172}} &
\mbox{\ding{172}} \\ 0 & 0  & -2\mbox{\ding{172}}-1 &
-2\mbox{\ding{172}}+3 \\ 0 & 2 & -\mbox{\ding{172}}+3 &
-\mbox{\ding{172}}+1 \end{array}  \right]
\]
 \[
% \hspace*{-7mm}
 \left[ \begin{array}{ccc|c}  1 & 0  & \mbox{\ding{172}} &
\mbox{\ding{172}} \\ 0 & 1  &
-2\mbox{\ding{172}}^2-\mbox{\ding{172}} &
-2\mbox{\ding{172}}^2+3\mbox{\ding{172}} \\ 0 & 2 &
-\mbox{\ding{172}}+3 & -\mbox{\ding{172}}+1 \end{array}  \right]
 \rightarrow
 \left[ \begin{array}{ccc|c}  1 & 0  & \mbox{\ding{172}} &
\mbox{\ding{172}} \\ 0 & 1  &
-2\mbox{\ding{172}}^2-\mbox{\ding{172}} &
-2\mbox{\ding{172}}^2+3\mbox{\ding{172}} \\ 0 & 0 &
4\mbox{\ding{172}}^2+\mbox{\ding{172}}+3 &
4\mbox{\ding{172}}^2-7\mbox{\ding{172}}+1 \end{array}  \right]
  \]
\[
\hspace*{-7mm}
 \left[ \begin{array}{ccc|c}  1 & 0  & \mbox{\ding{172}} &
\mbox{\ding{172}} \\ 0 & 1  &
-2\mbox{\ding{172}}^2-\mbox{\ding{172}} &
-2\mbox{\ding{172}}^2+3\mbox{\ding{172}} \\ 0 & 0 & 1 &
\frac{4\mbox{\ding{172}}^2-7\mbox{\ding{172}}+1}{4\mbox{\ding{172}}^2+\mbox{\ding{172}}+3}
\end{array}  \right]
 \rightarrow
 \left[ \begin{array}{ccc|c}  1 & 0  & 0 &
\frac{8\mbox{\ding{172}}^2+2\mbox{\ding{172}}}{4\mbox{\ding{172}}^2+\mbox{\ding{172}}+3}
\\ 0 & 1  & 0 &
\frac{-8\mbox{\ding{172}}^2+10\mbox{\ding{172}}}{4\mbox{\ding{172}}^2+\mbox{\ding{172}}+3}
\\ 0 & 0 & 1 &
\frac{4\mbox{\ding{172}}^2-7\mbox{\ding{172}}+1}{4\mbox{\ding{172}}^2+\mbox{\ding{172}}+3}
\end{array}  \right]
\]
It is easy to see that the finite parts of the numbers
\[
  \tilde{x}^*_1 = \frac{8\mbox{\ding{172}}^2+2\mbox{\ding{172}}}{4\mbox{\ding{172}}^2+\mbox{\ding{172}}+3}
= 2- \frac{6}{4\mbox{\ding{172}}^2+\mbox{\ding{172}}+3},
  \]
\[  \tilde{x}^*_2 = \frac{-8\mbox{\ding{172}}^2+10\mbox{\ding{172}}}{4\mbox{\ding{172}}^2+\mbox{\ding{172}}+3}
 = -2 +
\frac{12\mbox{\ding{172}}+6}{4\mbox{\ding{172}}^2+\mbox{\ding{172}}+3},
  \]
\[
      \tilde{x}^*_3 = \frac{4\mbox{\ding{172}}^2-7\mbox{\ding{172}}+1}{4\mbox{\ding{172}}^2+\mbox{\ding{172}}+3}
 = 1-
\frac{8\mbox{\ding{172}}+2}{4\mbox{\ding{172}}^2+\mbox{\ding{172}}+3},
   \]
coincide with the corresponding solution $ x^*_1= 2,$ $     x^*_2
= -2,$ and $    x^*_3  = 1$.

In this procedure we have introduced   the number
$\mbox{\ding{172}}^{-1}$ two times. As a result, we have obtained
expressions where the maximal power of grossone is equal to 2 and
there are rational expressions depending on grossone, as well. By
reasoning analogously to Example~\ref{e_calcolo_3},    when we
  execute divisions, in the obtained
results it is not necessary to store the parts of the type
$c_p\mbox{\ding{172}}^p, p < -2,$   because in any way the result
of their successive multiplication with the numbers of the type
$c_2\mbox{\ding{172}}^2$   will give finite exponents less than
zero. That is, numbers with these exponents will be infinitesimals
that are not interesting for us in this computational context.
Thus, by using the positional numeral system (\ref{3.12}),
(\ref{3.13}), we obtain
\[
\hspace*{-2mm}
 \left[ \begin{array}{ccc|c}  1 & 0  & \mbox{\ding{172}} &
\mbox{\ding{172}} \\ 0 & 1  &
-2\mbox{\ding{172}}^2-\mbox{\ding{172}} &
-2\mbox{\ding{172}}^2+3\mbox{\ding{172}} \\ 0 & 0 & 1 &
\frac{4\mbox{\ding{172}}^2-7\mbox{\ding{172}}+1}{4\mbox{\ding{172}}^2+\mbox{\ding{172}}+3}
\end{array}  \right]
 \rightarrow
 \left[ \begin{array}{ccc|c}  1 & 0  & \mbox{\ding{172}} &
\mbox{\ding{172}} \\ 0 & 1  &
-2\mbox{\ding{172}}^2-\mbox{\ding{172}} &
-2\mbox{\ding{172}}^2+3\mbox{\ding{172}} \\ 0 & 0 & 1 &
1\mbox{\ding{172}}^0\mbox{\small{-}}2\mbox{\ding{172}}^{-1}
\end{array}  \right].
\]
Note that the number
$1\mbox{\ding{172}}^0\mbox{\small{-}}2\mbox{\ding{172}}^{-1}$ does
not contain the part of the type $c_{-2}\mbox{\ding{172}}^{-2}$
because the coefficient $c_{-2}$ obtained after the executed
division is such that $c_{-2}=0$. Then we proceed as follows
\[
\hspace*{-2mm}
 \left[ \begin{array}{ccc|c}  1 & 0  & \mbox{\ding{172}} &
\mbox{\ding{172}} \\ 0 & 1  & 0 & -2 \\ 0 & 0 & 1 &
1\mbox{\ding{172}}^0\mbox{\small{-}}2\mbox{\ding{172}}^{-1}\end{array}
\right]
 \rightarrow
 \left[ \begin{array}{ccc|c}  1 & 0  & 0 &
2 \\ 0 & 1  & 0 & -2 \\ 0 & 0 & 1 &
1\mbox{\ding{172}}^0\mbox{\small{-}}2\mbox{\ding{172}}^{-1}
\end{array}  \right].
\]
The obtained solutions  $ x^*_1= 2$ and $x^*_2 = -2$ have been
obtained exactly without infinitesimal parts and $ x^*_3 = 1$ is
derived from the finite part of
$1\mbox{\ding{172}}^0\mbox{\small{-}}2\mbox{\ding{172}}^{-1}$.
\hfill $\Box$
\end{example}

We conclude this section by emphasizing that zero pivots in the
matrix are substituted dynamically by $\mbox{\ding{172}}^{-1}$.
Thus, the number of the introduced infinitesimals
$\mbox{\ding{172}}^{-1}$ depends on the number of zero pivots.

\section{New computational possibilities for mathematical\\
modelling} \label{s_calcolo4}

The computational capabilities of the Infinity Computer allow one
to construct new  and more powerful  mathematical models  able to
take into account infinite and infinitesimal changes of
parameters. In this section, the main attention is given to
infinitesimals that can increase the accuracy of models and
computations, in general. It is  shown that the introduced
infinitesimal numerals allow us to formalize the concept `point'
and to use it in practical calculations. Examples related to
computations of probabilities and areas (and volumes) of objects
having several parts of different dimensions are given.

\subsection{Numerical representations of points at an interval}
\label{s4.2}

We start by reminding   traditional definitions of the infinite
sequences and subsequences.  An \textit{infinite sequence}
$\{a_n\}, a_n \in A, n \in \mathbb{N},$ is a function having as
the domain the set of natural numbers, $\mathbb{N}$, and as the
codomain  a set $A$. A \textit{subsequence} is   a sequence from
which some of its elements have been removed.
\begin{theorem}
\label{t2} The number of elements of any infinite sequence is less
or equal to~\ding{172}.
\end{theorem}

\textit{Proof.}  It has been postulated that the set $\mathbb{N}$
has \ding{172} elements. Thus, due to the sequence definition
given above, any sequence having $\mathbb{N}$ as the domain  has
\ding{172} elements.

The notion of subsequence is introduced as a sequence from which
some of its elements have been removed. Thus, this definition
gives infinite sequences having the number of members less than
grossone.  \hfill $\Box$

It becomes appropriate now to define the \textit{complete
sequence} as an infinite sequence  containing \ding{172} elements.
For example, the sequence   of natural numbers  is complete, the
sequences of even  and odd natural numbers  are not complete. One
of the immediate consequences of the understanding of this result
is that any sequential process can have at maximum \ding{172}
elements and, due to Postulate 1, it depends on the chosen numeral
system which numbers among  \ding{172} members of the process we
can observe.

By using the introduced, more precise than the traditional one,
definition of   sequence, we can  calculate the number of points
of the interval $[0,1)$, of a line, and of the $N$-dimensional
space. To do this we need a definition of the term
`point'\index{point} and mathematical tools to indicate a point.
Since this concept is one of the most fundamental, it is very
difficult to find an adequate definition. If we accept (as is
usually done in modern Mathematics) that a \textit{point} $A$
belonging to the interval $[0,1)$ is determined by a numeral $x$,
$x  \in \mathbb{S},$ called \textit{coordinate  of the point A}
where $\mathbb{S}$ is a set of numerals,    then we can indicate
the point $A$ by its coordinate  $x$  and we are able to execute
the required calculations.

It is worthwhile to emphasize  that we have not postulated that
$x$ belongs to the   set, $\mathbb{R}$, of real numbers  as it is
usually done, because we can express coordinates only by numerals
and different choices of numeral systems lead to various sets of
numerals. This situation   is a direct consequence of Postulate~2
and is typical for natural sciences where it is well known that
instruments influence the result  of observations. It is similar
as to   work with a microscope: we decide the level of the
precision we need and obtain a result which is dependent on the
chosen level of accuracy. If we need a more precise or a more
rough answer, we change the lens of our microscope.

We should decide now which numerals we shall use to express
coordinates of the points. Different variants can be chosen
depending on the precision level we want to obtain. For example,
if the numbers $0 \le x < 1$ are expressed in the form
$\frac{p-1}{\mbox{\ding{172}}}, p  \in \mathbb{N}$, then the
smallest positive number we can distinguish is
$\frac{1}{\mbox{\ding{172}}}$. Therefore, the interval $[0,1)$
contains the  following \ding{172} points
\[
0, \,\,\,\, \frac{1}{\mbox{\ding{172}}}, \,\,\,\,
\frac{2}{\mbox{\ding{172}}}, \,\,\,\,  \ldots \,\,\,\,
\frac{\mbox{\ding{172}}-2}{\mbox{\ding{172}}}, \,\,\,\,
\frac{\mbox{\ding{172}}-1}{\mbox{\ding{172}}}.
\]
Then, due to Theorem~\ref{t2}  and the definition of sequence,
  \ding{172} intervals of the form $[a-1,a), a \in
\mathbb{N},$ can be distinguished at the ray  $x \ge 0$. Hence,
this ray contains $\mbox{\ding{172}}^{2}$ points and the whole
line consists of $2\mbox{\ding{172}}^{2}$ points.

If we need a higher precision, within each interval
\[
[(a-1)\frac{i-1}{\mbox{\ding{172}}},a\frac{i}{\mbox{\ding{172}}}),
\hspace{3mm} a,i \in \mathbb{N},
\]
we can distinguish again \ding{172} points and the number of
points within each interval $[a-1,a), a \in \mathbb{N},$ will
become  equal to $\mbox{\ding{172}}^{2}$. Consequently, the number
of the points of this kind  on the line will be equal to
$2\mbox{\ding{172}}^{3}$.

Continuing the analogy with the microscope, we can also decide to
change our microscope with a new one. In our terms this means to
change the numeral system with another one. For instance, instead
of the numerals considered above, we choose a positional numeral
system to calculate the number of points within the interval
$[0,1)$ expressed by numerals
 \beq
(.a_{-1} a_{-2}  \ldots a_{-(\mbox{\tiny\ding{172}}-1)}
a_{-\mbox{\tiny\ding{172}}})_b.
 \label{3.103calcolo}
       \eeq

\begin{theorem}
\label{t5calcolo} The number of elements   of the set
 of numerals (\ref{3.103calcolo}) is equal to
$ b^{\mbox{\tiny\ding{172}}}$.
\end{theorem}

\textit{Proof.}  The proof is obvious and is so omitted. \hfill
 $\Box$
\begin{corollary}
The  number of points expressed by numerals
 \beq
 (a_{\mbox{\tiny\ding{172}}-1}a_{\mbox{\tiny\ding{172}}-2} \ldots a_1 a_0 .
 a_{-1} a_{-2}  \ldots   a_{-(\mbox{\tiny\ding{172}}-1)}     a_{-\mbox{\tiny\ding{172}}})_b
   \label{3.103}
       \eeq
is equal to $b^{2\mbox{\tiny\ding{172}}}$.
\end{corollary}

\textit{Proof.}  The corollary is a straightforward consequence of
Theorem~\ref{t5calcolo}. \hfill
 $\Box$

\subsection{Applications in probability theory and calculating volumes}
\label{s4.6}

A formalization of the concept `point' introduced above allows us
to execute  more accurately computations having relations with
this concept. Very often in   scientific computing and engineering
it is required to construct mathematical models for
multi-dimensional objects. Usually this is done by partitioning
the modelled object in several parts having the same dimension and
each of the parts is modelled separately. Then additional efforts
are made in order to provide a correct functioning of a model
unifying the obtained sub-models and  describing the entire
object.

Another interesting applied area is linked to stochastic models
dealing with events having probability equal to zero. In this
subsection, we first show that the new approach allows us to
distinguish   the impossible event having the probability equal to
zero (i.e., $P(\varnothing)=0$) and events having an infinitesimal
probability. Then we show how infinitesimals can be used in
calculating volumes of objects consisting of parts having
different dimensions.

\begin{figure}[t]
  \begin{center}
    \epsfig{ figure = 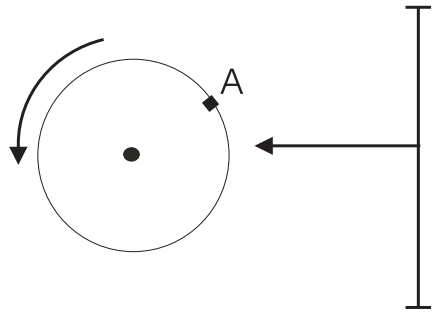, width = 2.2in, height = 1.6in,  silent = yes }
    \caption{What is the probability that the rotating disk stops in
such a way that the point $A$ will be exactly in front of the
arrow? }
 \label{Big_paper2}
  \end{center}
\end{figure}

 Let us consider the problem presented in
Fig.~\ref{Big_paper2} from the traditional point of view of
probability theory. We start to rotate a disk having radius $r$
with the point $A$ marked at its border and we would like to know
the  probability $P(E)$ of the following event $E$:  the disk
stops in such a way that the point $A$ will be exactly in front of
the arrow fixed at the wall.   Since the point $A$  is an entity
that has no extent it is calculated by considering the following
limit
\[
P(E) = \lim_{h \rightarrow 0}\frac{h}{2\pi r}=0.
\]
where $h$ is an arc of the circumference containing $A$ and $2\pi
r$ is its length.

However,   the point $A$ can stop in front of the arrow, i.e.,
this event is not impossible and its probability should be
strictly greater than zero, i.e., $P(E)>0$. The new approach
allows us to calculate   this probability numerically.

First of all, in order to state the experiment more rigorously, it
is necessary to choose a numeral system to express the points on
the circumference. This choice will fix the number of points, $K$,
that we are able to distinguish  on the circumference. Definition
of the notion \textit{point} allows us to define elementary events
in our experiment as follows:  the disk has stopped and the arrow
indicates a point. As a consequence, we obtain that the number,
$N(\Omega)$, of all possible elementary events, $e_i$,  in our
experiment is equal to $K$ where
$\Omega=\cup_{i=1}^{N(\Omega)}e_i$ is the sample space of our
experiment. If our disk is well balanced, all elementary events
are equiprobable and, therefore, have the same probability equal
to $\frac{1}{N(\Omega)}$. Thus, we can calculate $P(E)$ directly
by subdividing the number, $N(E)$, of favorable elementary events
by the number, $K=N(\Omega)$, of all possible events.

For example, if we use numerals of the type
$\frac{i}{\mbox{\ding{172}}}, i \in \mathbb{N},$ then
$K=\mbox{\ding{172}}$. The number $N(E)$ depends on our decision
about how many numerals we want to use to represent the point $A$.
If we decide that the point~$A$ on the circumference is
represented by $m$ numerals  we obtain
\[
P(E) = \frac{N(E)}{N(\Omega)}= \frac{m}{K} =
\frac{m}{\mbox{\ding{172}}} > 0.
\]
where the number $\frac{m}{\mbox{\ding{172}}}$ is infinitesimal if
$m$ is finite. Note that this representation is very interesting
also from the point of view of distinguishing   the notions
`point' and `arc'. When $m$ is finite than we deal with a point,
when $m$ is infinite we deal with an arc.

In the case we need a higher accuracy, we can choose, for
instance, numerals of the type $ i \mbox{\ding{172}}^{-2}, 1 \le i
\le \mbox{\ding{172}}^{2},$ for expressing points at the disk.
Then it follows $K=\mbox{\ding{172}}^2$ and, as a result, we
obtain $P(E) = m \mbox{\ding{172}}^{-2}  > 0$.

This example with the rotating disk, of course, is a particular
instance of the general
 situation taking place in the traditional probability
theory where the probability that  a continuous random variable
$X$ attains a given value $a$ is zero, i.e., $P(X=a)=0$. While for
a discrete random variable one could say that an event with
probability zero is impossible, this can not be said in the case
of a continuous random variable. As we have shown by the example
above, in our approach this situation does not take place because
this probability can be expressed by infinitesimals. As a
consequence, probabilities of such events can be computed and used
in numerical models describing the real world  (see \cite{chaos}
for a detailed discussion on the modelling continuity by
infinitesimals in the framework of the  approach using grossone).

Moreover, the obtained probabilities are not absolute,
\textit{they depend on the accuracy chosen for the mathematical
model describing the experiment}. There is again a straight
analogy with Physics  where it is not possible to obtain results
that have a precision higher than the accuracy of the measurement
of the  data. We also cannot obtained  a precision that is higher
than the precision of numerals used in the mathematical model.

 \begin{figure}[t]
  \begin{center}
    \epsfig{ figure = 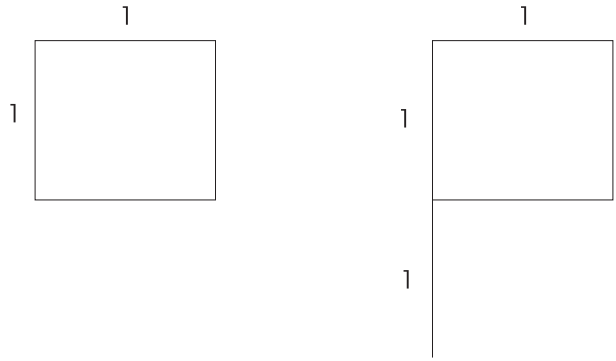, width = 8cm,height =47mm,  silent = yes }
    \caption{It is possible to calculate and to distinguish areas of these two objects }
 \label{Big_paper9}
  \end{center}
\end{figure}

\begin{figure}[t]
  \begin{center}
    \epsfig{ figure = 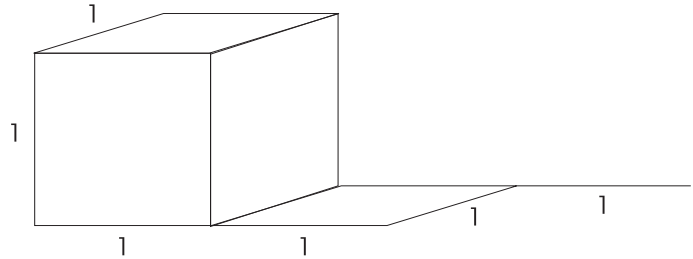, width = 9cm,height =34mm,  silent = yes }
    \caption{New possibilities for calculating volumes of objects}
 \label{Big_paper10}
  \end{center}
\end{figure}

Let us now consider two examples showing that the new approach
allows us to calculate areas and volumes of a more general class
of objects than the traditional one. In Fig.~\ref{Big_paper9}
  two figures are shown. The traditional approach tells us that both
of them have area equal to one. In the new approach,    if we use
numerals of the type $i\mbox{\ding{172}}^{-1}, i \in \mathbb{N},$
to express points within a unit interval, then the unit interval
consists of \ding{172} points and in the plane each point has the
infinitesimal area
$\mbox{\ding{172}}^{-1}\cdot\mbox{\ding{172}}^{-1}=\mbox{\ding{172}}^{-2}$.
As a consequence, this value will be our accuracy in calculating
areas in this example. Suppose now that the vertical line added to
the square at the right figure in Fig.~\ref{Big_paper9} has the
width equal to one point. Then we are able to calculate the area,
$S_2$, of the right figure and it will be possible to distinguish
it from the area, $S_1$, of the square on the left
 \[
 S_1 = 1 \cdot 1= 1, \hspace{1cm}  S_2 = 1\cdot 1 + 1 \cdot
\mbox{\ding{172}}^{-1}=1\mbox{\ding{172}}^{0}1\mbox{\ding{172}}^{-1}.
\]
If the added vertical line has the width equal to three points
then it follows
 \[
  S_2 = 1\cdot 1 + 3 \cdot
\mbox{\ding{172}}^{-1}=1\mbox{\ding{172}}^{0}3\mbox{\ding{172}}^{-1}.
\]
The  volume of the figure shown in Fig.~\ref{Big_paper10} is
calculated analogously:
 \[
\hspace{1cm} V = 1 \cdot 1 \cdot 1 + 1 \cdot 1 \cdot
\mbox{\ding{172}}^{-1} + 1\cdot  \mbox{\ding{172}}^{-1} \cdot
\mbox{\ding{172}}^{-1}
=1\mbox{\ding{172}}^{0}1\mbox{\ding{172}}^{-1}1
\mbox{\ding{172}}^{-2}.
\]
If the  accuracy of the considered numerals of the type
$i\mbox{\ding{172}}^{-1}, i \in \mathbb{N},$ is not sufficient, we
can increase it by using, for instance, numerals of the type
$i\mbox{\ding{172}}^{-2}, 1 \le i \le \mbox{\ding{172}}^{2}$. Then
the unit interval consists of $\mbox{\ding{172}}^{2}$ points and
at the plane each point has the infinitesimal area
$\mbox{\ding{172}}^{-2}\cdot\mbox{\ding{172}}^{-2}=\mbox{\ding{172}}^{-4}$.
As a result, by a complete analogy to the previous case we obtain
for lines having the width, for instance, equal to five points in
all three dimensions that
 \[
  S_2 = 1\cdot 1 + 5 \cdot
\mbox{\ding{172}}^{-2}=1\mbox{\ding{172}}^{0}5\mbox{\ding{172}}^{-2},
\]
\[
V = 1 \cdot 1 \cdot 1 + 1 \cdot 1 \cdot
5\cdot\mbox{\ding{172}}^{-2} + 1\cdot 5\cdot\mbox{\ding{172}}^{-2}
\cdot 5\cdot\mbox{\ding{172}}^{-2}
=1\mbox{\ding{172}}^{0}5\mbox{\ding{172}}^{-2}25
\mbox{\ding{172}}^{-4}.
\]

Finally, we conclude the paper by a remark that thanks to the new
approach it becomes possible to measure fractal objects at
infinity. A detailed discussion on measuring fractals by the
introduced infinite and infinitesimal numbers can be found in
\cite{chaos,spirals}.

%\markboth{Bibliography}{Bibliography}
\bibliographystyle{amsplain}
\bibliography{XBib_Calcolo}

\providecommand{\bysame}{\leavevmode\hbox to3em{\hrulefill}\thinspace}
\providecommand{\MR}{\relax\ifhmode\unskip\space\fi MR }
\providecommand{\MRhref}[2]{%
  \href{http://www.ams.org/mathscinet-getitem?mr=#1}{#2}
}
\providecommand{\href}[2]{#2}
\begin{thebibliography}{10}

\bibitem{Benci}
V.~Benci and M.~{Di~Nasso}, \emph{Numerosities of labeled sets: a new way of
  counting}, Advances in Mathematics \textbf{173} (2003), 50--67.

\bibitem{Cantor}
G.~Cantor, \emph{Contributions to the founding of the theory of transfinite
  numbers}, Dover Publications, New York, 1955.

\bibitem{Conway}
J.H. Conway and R.K. Guy, \emph{The book of numbers}, Springer-Verlag, New
  York, 1996.

\bibitem{Gordon}
P.~Gordon, \emph{Numerical cognition without words: {E}vidence from
  {A}mazonia}, Science \textbf{306} (2004), no.~15 October, 496--499.

\bibitem{Knopp}
K.~Knopp, \emph{Theory and application of infinite series}, Dover Publications,
  New York, 1990.

\bibitem{Loeb}
P.A. Loeb and M.P.H. Wolff, \emph{Nonstandard analysis for the working
  mathematician}, Kluwer Academic Publishers, Dordrecht, 2000.

\bibitem{Robinson}
A.~Robinson, \emph{Non-standard analysis}, Princeton Univ. Press, Princeton,
  1996.

\bibitem{Sergeyev}
{Ya.D.} Sergeyev, \emph{Arithmetic of infinity}, Edizioni Orizzonti
  Meridionali, CS, 2003.

\bibitem{Sergeyev_patent}
Ya.D. Sergeyev, \emph{Computer system for storing infinite, infinitesimal, and
  fi\-ni\-te quan\-ti\-ties and executing arithmetical operations with them},
  patent application 08.03.04, 2004.

\bibitem{www}
\bysame, \emph{http://www.theinfinitycomputer.com}, 2004.

\bibitem{Poland}
\bysame, \emph{Mathematical foundations of the {I}nfinity {C}omputer}, Annales
  {UMCS} {I}nformatica {AI} \textbf{4} (2006), 20--33.

\bibitem{Mathesis}
\bysame, \emph{Misuriamo l'infinito}, Periodico di Matematiche \textbf{6(2)}
  (2006), 11--26.

\bibitem{chaos}
\bysame, \emph{Blinking fractals and their quantitative analysis using infinite
  and infinitesimal numbers}, Chaos, Solitons $\&$ Fractals \textbf{33(1)}
  (2007), 50--75.

\bibitem{spirals}
\bysame, \emph{Measuring fractals by infinite and infinitesimal numbers},
  Mathematical Methods, Physical Methods $\&$ Simulation Science and Technology
  \textbf{1(1)} (2008), 217--237.

\bibitem{informatica}
\bysame, \emph{A new applied approach for executing computations with infinite
  and infinitesimal quantities}, Informatica (2008), (to appear).

\bibitem{Strongin_Sergeyev}
R.G. Strongin and Ya.D. Sergeyev, \emph{Global optimization and non-convex
  constraints: Sequential and parallel algorithms}, Kluwer Academic Publishers,
  Dordrecht, 2000.

\bibitem{Strongin&Sergeyev(2003)}
\bysame, \emph{Global optimization: {F}ractal approach and non-redundant
  parallelism}, J. Global Optimization \textbf{27(1)} (2003), 25--50.

\bibitem{Walster_1}
G.W. Walster, \emph{Compiler support of interval arithmetic with inline code
  generation and nonstop exception handling}, Tech. Report, Sun Microsystems,
  2000.

\end{thebibliography}
%\end{article}
\end{document}